\documentclass[a4paper,12pt]{amsart}
\usepackage{amsfonts}
\usepackage{amssymb}
\usepackage{ifthen}
\usepackage{graphicx}
\usepackage{color}
\nonstopmode \numberwithin{equation}{section}
\newtheorem{thm}{Theorem}
\newtheorem{lem}{Lemma}
\newtheorem{cor}{Corollary}
\newtheorem{prop}{Proposition}

\newtheorem{conj}{Conjecture}

\theoremstyle{definition}
\newtheorem{defn}{Definition}
\newtheorem{example}{Example}

\newtheorem{ques}{Question}
\newtheorem{rem}{Remark}

\newtheorem{rems}{Remarks}

\newcounter {own}
\def\theown {\thesection  .\arabic{own}}

\newenvironment{pf}[1][]{%
 \vskip 3mm
 \noindent
 \ifthenelse{\equal{#1}{}}%
  {{\slshape Proof. }}%
  {{\slshape #1.} }%
 }%
{\qed\bigskip}

\newcounter{alphabet}
\newcounter{tmp}
\newenvironment{Thm}[1][]{\refstepcounter{alphabet}%
\medskip%
\noindent%
{\bf Theorem \Alph{alphabet}}%
\ifthenelse{\equal{#1}{}}{}{ (#1)}%
{\bf .} \itshape}{\vskip 8pt}

\makeatletter
\newcommand{\Ref}[1]{\@ifundefined{r@#1}{}{\setcounter{tmp}{\ref{#1}}\Alph{tmp}}}
\makeatother

\newcommand{\N}{{\mathbb N}}

\def\be{\begin{equation}}
\def\ee{\end{equation}}

\newcommand{\bee}{\begin{enumerate}}
\newcommand{\eee}{\end{enumerate}}

\newcommand{\blem}{\begin{lem}}
\newcommand{\elem}{\end{lem}}
\newcommand{\bthm}{\begin{thm}}
\newcommand{\ethm}{\end{thm}}
\newcommand{\bcor}{\begin{cor}}
\newcommand{\ecor}{\end{cor}}
\newcommand{\beg}{\begin{example}}
\newcommand{\eeg}{\end{example}}
\newcommand{\begs}{\begin{examples}}
\newcommand{\eegs}{\end{examples}}
\newcommand{\bdefn}{\begin{defn}}
\newcommand{\edefn}{\end{defn}}
\newcommand{\bprob}{\begin{prob}}
\newcommand{\eprob}{\end{prob}}
\newcommand{\bei}{\begin{itemize}}
\newcommand{\eei}{\end{itemize}}
\newcommand{\bqn}{\begin{ques}}
\newcommand{\eqn}{\end{ques}}
\newcommand{\bcon}{\begin{conj}}
\newcommand{\econ}{\end{conj}}
\newcommand{\bcons}{\begin{conjs}}
\newcommand{\econs}{\end{conjs}}
\newcommand{\bprop}{\begin{prop}}
\newcommand{\eprop}{\end{prop}}
\newcommand{\brem}{\begin{rem}}
\newcommand{\erem}{\end{rem}}
\newcommand{\brems}{\begin{rems}}
\newcommand{\erems}{\end{rems}}
\newcommand{\bo}{\begin{obser}}
\newcommand{\eo}{\end{obser}}
\newcommand{\bos}{\begin{obsers}}
\newcommand{\eos}{\end{obsers}}
\newcommand{\bpf}{\begin{pf}}
\newcommand{\epf}{\end{pf}}
\newcommand{\ba}{\begin{array}}
\newcommand{\ea}{\end{array}}
\newcommand{\beq}{\begin{eqnarray}}
\newcommand{\beqq}{\begin{eqnarray*}}
\newcommand{\eeq}{\end{eqnarray}}
\newcommand{\eeqq}{\end{eqnarray*}}

\newcommand{\ds}{\displaystyle}

\newcounter{minutes}
\divide\time by 60
\newcounter{hours}
\multiply\time by 60 \addtocounter{minutes}{-\time}

\begin{document}
\bibliographystyle{amsplain}
\title[Multiplication operators between discrete Hardy spaces on rooted trees]
{Multiplication operators between discrete Hardy spaces on rooted trees}

%

\author{P. Muthukumar}
\address{P. Muthukumar, Indian Statistical Institute,
Statistics and Mathematics Unit, 8th Mile, Mysore Road,
Bangalore, 560 059, India.}
\email{pmuthumaths@gmail.com}

\author{P. Shankar}
\address{P. Shankar, Indian Statistical Institute,
Statistics and Mathematics Unit, 8th Mile, Mysore Road,
Bangalore, 560 059, India.}
\email{shankarsupy@gmail.com }




\begin{abstract}
Muthukumar and Ponnusamy \cite{MP-Tp-spaces} studied the multiplication operators on $\mathbb{T}_p$ spaces. In this article, we mainly consider multiplication operators between $\mathbb{T}_p$ and $\mathbb{T}_q$ ($p\neq q$). In particular, we characterize bounded and compact multiplication operators from  $\mathbb{T}_{p}$ to $\mathbb{T}_{q}$. For $p\neq q$, we prove that there are no invertible multiplication operators from $\mathbb{T}_{p}$ to $\mathbb{T}_{q}$ and also there are no isometric multiplication operators from $\mathbb{T}_{p}$ to $\mathbb{T}_{q}$. Finally,  we discuss about fixed points of a multiplication operator on $\mathbb{T}_{p}$.
\end{abstract}
\thanks{
File:~\jobname .tex,
          printed: \number\day-\number\month-\number\year,
          \thehours.\ifnum\theminutes<10{0}\fi\theminutes
}
\maketitle
\pagestyle{myheadings}
\markboth{P. Muthukumar and P. Shankar}{Multiplication operators between discrete Hardy spaces on rooted trees}

\section{Introduction}
Let $X$ be a normed space of functions defined on a set $\Omega$. For a complex-valued map $\psi$  on $\Omega$, the inducing multiplication operator is denoted by $M_\psi$ and defined by
$$M_\psi f = \psi f \mbox{~for every~}f \in X.$$

The study of multiplication operators acts as bridge between operator theory and function theory. Multiplication operators on various analytic function spaces on the open unit disk have been studied for many years, and the literature is extensive. Recent years, considerable work has been done on multiplication operators on the spaces of functions defined on trees. See for example \cite{Colonna-MO-5,Colonna-MO-3,Colonna-MO-2,Colonna-MO-6,Allen-MO-7,Colonna-MO-1,Colonna-MO-4}.

A discrete analogue ($\mathbb{T}_{p}$) of Hardy spaces on homogeneous rooted trees were defined in \cite{MP-Tp-spaces}. In the same article, multiplication operators on $\mathbb{T}_{p}$ spaces were studied. Study of composition operators on $\mathbb{T}_{p}$ spaces were intiated and developed in \cite{CO-Tp-spaces} and \cite{CO-Tp-spaces2}, respectively.
In this article, we deal with the study of multiplication operators between various $\mathbb{T}_{p}$ spaces. Also, we consider general rooted trees instead of homogeneous rooted trees.

The paper is organized as follows. We refer to Section \ref{Msha1-Sec2} for basic definitions and preliminaries about $\mathbb{T}_{p}$ spaces. In Section \ref{Msha1-Sec3}, we give characterization of bounded multiplication operators from $\mathbb{T}_{p}$   to $\mathbb{T}_{q}$. Also, we compute their operator norms for each case. We investigate about isometric multiplication operators from $\mathbb{T}_{p}$ to $\mathbb{T}_{q}$  in Section \ref{Msha1-Sec4}. In Section \ref{Msha1-Sec5},
we provide necessary and sufficient conditions for a multiplication operator
to be a compact operator. In Section \ref{Msha1-Sec6}, we prove that there is no invertible multiplication operators from $\mathbb{T}_{p}$ to $\mathbb{T}_{q}$ for $p\neq q$. Finally, in Section \ref{Msha1-Sec7}, we discuss about fixed points of a multiplication operator on $\mathbb{T}_{p}$.

\section{Preliminaries}\label{Msha1-Sec2}

Let us fix some notation for throughout the paper. Let $T$ be an infinite rooted tree with root $\textsl{o}$. For each vertex $v\in T$, the number of edges in the unique path connecting the root $\textsl{o}$ and the vertex $v$, is denoted  by $|v|$.
For $n\in \N_0$, let $D_n$ denote the set of all vertices $v\in T$ with $|v| = n$ and
denote the number of elements in $D_n$ by $c_n$. Further, $\mathbb{C}$ denotes the set of all complex numbers, $\mathbb{N}=\{1,2,\ldots \}$  and $\mathbb{N}_0=\mathbb{N}\cup \{0\}$. By a function on $T$, we mean a complex-valued function defined on the vertex set of $T$.

For every $p\in(0,\infty]$, the discrete analogue of the generalized Hardy space $\mathbb{T}_{p}$ (see \cite{MP-Tp-spaces}), is defined by
$$\mathbb{T}_{p}:=\{f\colon T \to\mathbb{C} \, \big |\, \|f\|_{p}:= \sup\limits_{n\in \mathbb{N}_{0}} M_{p}(n,f)<\infty\},
$$
where $M_{\infty}(n,f):= \max\limits_{|v|=n } |f(v)|$ and for $p<\infty$,
$$
M_{p}(n,f):= \ds \left (\frac{1}{c_n}\sum\limits_{|v|=n}|f(v)|^{p} \right )^{\frac{1}{p}},
$$
for every $n\in \mathbb{N}_0$.

Similarly, the discrete analogue of the generalized little Hardy space, denoted by $\mathbb{T}_{p,0}$, is defined by
$$\mathbb{T}_{p,0}:=\{f\in\mathbb{T}_{p} :\, \lim\limits_{n\rightarrow\infty} M_{p}(n,f)=0\},
$$
for every $p\in(0,\infty]$.

Two vertices $u,v\in T$ are said to be neighbours if there is an edge connecting them. If every
vertex of $T$ has $k$ neighbours, then $T$ is said to be $k$-homogeneous tree. The following results in this section were proved  in \cite{MP-Tp-spaces} for homogeneous rooted trees. It is trivial to see that the following  results are valid for general rooted trees and proofs are similar to that of homogeneous rooted trees case.

\begin{Thm}\label{lem:bound} \emph{(Growth Estimate)}
{\rm (\cite[Lemma 3.12]{MP-Tp-spaces})}
Let $T$ be an infinite rooted tree and $0<p<\infty$.
If $f$ is an element of $\mathbb{T}_{p}$ or $\mathbb{T}_{p,0}$, then we have
$$ |f(v)|\leq \ds{(c_{|v|})^{\frac{1}{p}}\, \|f\|_p} ~\mbox{ for all }~ v\in T.
$$
\end{Thm}

\begin{Thm}\label{thm:banachp}
{\rm (\cite[Theorems 3.1 and 3.5]{MP-Tp-spaces})}
For $1\leq p\leq\infty$, $\|.\|_{p}$ induces a Banach space structure on the spaces $\mathbb{T}_{p}$ and
$\mathbb{T}_{p,0}$.
\end{Thm}

\begin{Thm}
{\rm (\cite[Theorems 3.10 and 3.11]{MP-Tp-spaces})}
For $0<p\leq\infty$, the space $\mathbb{T}_{p}$ is not  separable, whereas $\mathbb{T}_{p,0}$ is a
separable space as the span of $\{\chi_v :\, v\in T\}$ is
dense in $\mathbb{T}_{p,0}$.
\end{Thm}

\section{bounded Multiplication Operators}\label{Msha1-Sec3}

Before we get into discussion on multiplication operators, let us observe from Lemma \ref{inequ}, the inclusion properties of  $\mathbb{T}_p$ spaces, namely:
\bee
\item If $\sup_{n\in \mathbb{N}_0} c_n< \infty$, then $\mathbb{T}_p=\mathbb{T}_\infty$ for all $p>0$.
\item If $\sup_{n\in \mathbb{N}_0} c_n= \infty$, then $\mathbb{T}_q\varsubsetneq\mathbb{T}_p$ for $0<p<q$.
\eee

\blem\label{inequ}
Let $T$ be an infinite rooted tree and $0<p<q<\infty$. Then, for $f:T\rightarrow \mathbb{C}$ and for each $n\in \mathbb{N}_0$, we have
$$
\ds\frac{M_\infty(n,f)}{{c_n}^\frac{1}{p}}\leq  \frac{1}{{c_n}^{\frac{1}{p}-\frac{1}{q}}} \,M_q(n,f)\leq M_p(n,f)\leq M_q(n,f)\leq M_\infty (n,f).
$$

\elem

\bpf
Let $f:T\rightarrow \mathbb{C}$ be an arbitrary function.
Fix $n\in \mathbb{N}_0$. If $|f(v)|\leq M \mbox{~for all~} v\in D_n$ for some $M>0$, then
$$ M_q(n,f)=\ds\left(\frac{1}{c_n}\sum\limits_{|v|=n}|f(v)|^q \right)^{\frac{1}{q}}\leq M.$$ Thus,

\begin{equation}
M_q(n,f)\leq \max\limits_{|v|=n } |f(v)|= M_\infty (n,f).
\end{equation}

By  \cite[Lemma 3.6]{MP-Tp-spaces}, we have for $0<p<q<\infty$,

\begin{equation}
M_p(n,f)\leq M_q(n,f).
\end{equation}

For $0<p<q<\infty$, one has
$$\ds\left(\sum\limits_{|v|=n}|f(v)|^q\right)^\frac{1}{q}\leq \left(\sum\limits_{|v|=n}|f(v)|^p\right)^\frac{1}{p},$$ which gives

\begin{equation}
\frac{1}{{c_n}^{\frac{1}{p}-\frac{1}{q}}}M_q(n,f)\leq M_p(n,f).
\end{equation}

Finally, for each $v\in D_n$, note that
$$\ds\frac{|f(v)|^q}{{c_n}}\leq \frac{1}{{c_n}}\sum_{|w|=n}|f(w)|^q $$
and thus
$$\ds\frac{|f(v)|}{{c_n}^\frac{1}{q}}\leq M_q(n,f) ~\mbox{for all}~ v\in D_n.$$

In particular,

\begin{equation}\label{eqn4}
\ds\frac{M_\infty(n,f)}{{c_n}^\frac{1}{q}}\leq M_q(n,f).
\end{equation}

From above inequalities, for each $n\in \mathbb{N}_0$, we have
$$
\ds\frac{M_\infty(n,f)}{{c_n}^\frac{1}{p}}\leq  \frac{1}{{c_n}^{\frac{1}{p}-\frac{1}{q}}} \,M_q(n,f)\leq M_p(n,f)\leq M_q(n,f)\leq M_\infty (n,f).
$$
This completes the proof.
\epf

\bcor
Let $T$ be a rooted tree. If $\sup c_n< \infty$ (for example, $2$-homogeneous trees), then all $\mathbb{T}_p$ spaces are equal as a set. That is,
$$
\mathbb{T}_p=\mathbb{T}_\infty=\left\lbrace f:T\rightarrow \mathbb{C}: \sup\limits_{v\in T} |f(v)|<\infty \right\rbrace \mbox{~for all~} ~ 0<p\leq \infty.
$$
\ecor

\bcor
If $\{c_n\}$ is unbounded, then $\mathbb{T}_q$ is a proper subset of $\mathbb{T}_p$ for $0<p<q$.
\ecor

\bpf
By lemma \ref{inequ},
$\mathbb{T}_q$ is a subset of $\mathbb{T}_p$ for $0<p<q$. It is easy to see that this inclusion is proper.
For example,  choose a sequence of vertices $\{v_n\}$ such that $|v_n|=n$ for all $n\in \mathbb{N}$. Fix $r>0$ such that $p<r<q$. Consider the function $f$
defined by
$$f(v)= \left \{\begin{array}{cl}
\ds(c_n)^{\frac{1}{r}} & \mbox{ if $v=v_n$ for some $n \in \mathbb{N}$ }\\
0 &\mbox{ elsewhere.}
\end{array}
\right .
$$
Then, $f\in \mathbb{T}_{p}$ but not in $\mathbb{T}_{q}$.
\epf


 In this section, we will consider the bounded multiplication operators between various $\mathbb{T}_{p}$ spaces. We will look into this case by case.

\begin{Thm}\label{bdd} {\rm (\cite[Theorem 4.3]{MP-Tp-spaces})}
Let $T$ be a rooted tree, $\psi$ be a complex-valued function on $T$ and $0< p \leq\infty$.
Then
 $M_\psi:\mathbb{T}_{p}\, (\mbox{or} \,\mathbb{T}_{p,0}) \rightarrow \mathbb{T}_{p} \,(\mbox{or} \,\mathbb{T}_{p,0})$ is a bounded operator if and only if  $\psi$ is a bounded function on $T$, i.e., $\psi \in \mathbb{T}_{\infty}$. Further, we also have $\|M_\psi\|=\|\psi\|_\infty$.
\end{Thm}

\bthm\label{thm1}
Let $\psi$ be a complex-valued function on $T$ and let $0<q<p<\infty$. Then the following are equivalent:
\begin{enumerate}
\item[(a)] The multiplication operator $M_{\psi}:\mathbb{T}_p\rightarrow \mathbb{T}_q$ is bounded.

\item[(b)] The multiplication operator $M_{\psi}:\mathbb{T}_{p,0}\rightarrow \mathbb{T}_{q,0}$ is bounded.

\item[(c)] $\ds{\psi \in \mathbb{T}_{\frac{pq}{p-q}}}$.
\end{enumerate}
Moreover,
$\|M_\psi\|=\|\psi\|_{\frac{pq}{p-q}}. $
\ethm

\bpf
Let $s=\frac{p}{q}>1$ and choose $t=\frac{p}{p-q}>1$, so that $\frac{1}{s}+\frac{1}{t}=1$.
For $n\in\N_0$,
\beqq
M_q^q(n,M_\psi f)&= & \frac{1}{c_n}\sum_{|v|=n} |\psi(v)|^q \,|f(v)|^q\\
                 &\leq & \frac{1}{c_n}\left (\sum\limits_{|v|=n}|\psi(v)|^{qt}\right)^{\frac{1}{t}} \left (\sum\limits_{|v|=n}|f(v)|^{qs}\right)^{\frac{1}{s}}\\
                 &= & \left (\frac{1}{c_n}\sum\limits_{|v|=n}|\psi(v)|^{qt}\right)^{\frac{1}{t}} \left (\frac{1}{c_n}\sum\limits_{|v|=n}|f(v)|^{p}\right)^{\frac{q}{p}}\\
		&=& M_{qt}^q(n,\psi)\, M_p^q(n,f).
\eeqq
Thus, we have
$$M_q(n,M_\psi f)\leq \ds M_{\frac{pq}{p-q}}(n,\psi)\, M_p(n,f) \mbox{~for every~} n\in \N_0.$$
This gives that, whenever $\ds{\psi \in \mathbb{T}_{\frac{pq}{p-q}}}$, $M_{\psi}$ maps $\mathbb{T}_{p}$ to $\mathbb{T}_{q}$ and $\mathbb{T}_{p,0}$ to $\mathbb{T}_{q,0}$ with
$$
\|M_\psi\|\leq \|\psi\|_{\frac{pq}{p-q}}.
$$
This proves the implications $(c)\Rightarrow (a)$ and  $(c)\Rightarrow (b)$.

For the converse part, for each $n\in \N_0$,
define $f_n$ on $T$ by
$$ f_n(v)=\left\{
\begin{array}{cl}
|\psi(v)|^{\frac{q}{p-q}}& \mbox{ if } v\in D_n\\
0 & \mbox{ otherwise}.
\end{array}
\right.
$$
Then,
\beqq
M_{\frac{pq}{p-q}}(n,\psi)\,M_p(n,f_n)&=& \left (\frac{1}{c_n}\sum\limits_{|v|=n}|\psi(v)|^{\frac{pq}{p-q}}   \right)^{\frac{p-q}{pq}} \left (\frac{1}{c_n}\sum\limits_{|v|=n}|\psi(v)|^{\frac{pq}{p-q}}\right)^{\frac{1}{p}}\\
                                      &=&\left (\frac{1}{c_n}\sum\limits_{|v|=n}|\psi(v)|^{\frac{pq}{p-q}}   \right)^{\frac{1}{q}}\\
                                          &=&\left (\frac{1}{c_n}\sum\limits_{|v|=n}\left(|\psi(v)|\,|\psi(v)|^{\frac{q}{p-q}}\right)^q   \right)^{\frac{1}{q}}\\
                                          &=&\left (\frac{1}{c_n}\sum\limits_{|v|=n}|\psi(v)|^{q}\,|f_n(v)|^q   \right)^{\frac{1}{q}}\\
          &=&M_q(n,M_\psi f_n).
\eeqq
Since $f_n(v)=0$ for $v\notin D_n$,  $M_p(n,f_n)=\|f_n\|_p$ and $M_q(n,M_\psi f_n)=\|\psi f_n\|_q.$  Thus, for each $n\in \N_0$, we have
$$\|\psi f_n\|_q=M_{\frac{pq}{p-q}}(n,\psi) \|f_n\|_p.$$
Thus,
$$\|M_\psi \|=\sup\left\lbrace\frac{\|\psi f\|_q}{\|f\|_p} :f\neq 0 \right\rbrace\geq \sup\limits_{n\in \N_0} M_{\frac{pq}{p-q}}(n,\psi) =\|\psi\|_{\frac{pq}{p-q}}.$$
Since each $f_n\in \mathbb{T}_{p,0}$ and as well as  in $\mathbb{T}_{p}$, we obtain the implications
$(a)\Rightarrow (c)$ and $(b)\Rightarrow (c)$. Moveover, it also gives that
$$\|M_\psi\|\geq \|\psi\|_{\frac{pq}{p-q}}. $$
It completes the proof.
\epf

\bthm\label{bdd2}
Let $\psi$ be a function from $T$ to $\mathbb{C}$ and let $0<p<\infty$. Then the following are equivalent:
\begin{enumerate}
\item[(a)] The multiplication operator $M_{\psi}:\mathbb{T}_\infty\rightarrow \mathbb{T}_p$ is bounded.

\item[(b)] The multiplication operator $M_{\psi}:\mathbb{T}_{\infty,0}\rightarrow \mathbb{T}_{p,0}$ is bounded.

\item[(c)] $\psi \in \mathbb{T}_p$.
\end{enumerate}
In this case,  $$\|M_\psi \|=\|\psi \|_p.$$
\ethm

\bpf
Since the proof of
$(c)\Leftrightarrow (b)$ is similar to that of $(c)\Leftrightarrow (a)$, it is enough to prove that $(c)\Leftrightarrow (a)$. To prove this, we see that for each $n\in \N_0$,
\beqq
M_p(n,\psi f)&=&\left(\frac{1}{c_n}\sum_{|v|=n}|\psi(v)|^p\,|f(v)|^p  \right)^\frac{1}{p}
             \leq  M_p(n,\psi)\,M_\infty (n,f).
\eeqq
By taking supremum over $n\in \N_0$ on both sides, we get that
$$ \|\psi f\|_p\leq \|\psi\|_p\,\|f\|_\infty.$$
It tells that, if $\psi \in \mathbb{T}_p$ then $M_{\psi}:\mathbb{T}_\infty\rightarrow \mathbb{T}_p$ is a bounded operator with $\|M_\psi \|\leq \|\psi \|_p $.

Let us prove the converse part now. Assume that $M_{\psi}:\mathbb{T}_\infty\rightarrow \mathbb{T}_p$ is bounded operator.

Fix $n\in \N_0$.
Define $f$ on $T$ by
$$ f(v)=\left\{
\begin{array}{cl}
1& \mbox{ if } v\in D_n\\
0 & \mbox{ otherwise},
\end{array}
\right.
$$
so that $\|f \|_\infty=1$ and
$$M_p(n,\psi f)=\left(\frac{1}{c_n}\sum_{|v|=n}|\psi(v)|^p|f(v)|^p\right)^\frac{1}{p}=\left(\frac{1}{c_n}\sum_{|v|=n}|\psi(v)|^p \right)^\frac{1}{p}=M_p(n,\psi).$$

Thus, $M_p(n,\psi) \leq \|M_\psi \|$. Since $n\in \N_0$ was arbitrary, we get

$$\|M_\psi \|\geq \sup_{n\in \mathbb{N}_0} M_p(n,\psi)=\|\psi \|_p.$$
Therefore, $\psi$ should be in $\mathbb{T}_p$.
This completes the proof.
\epf


\bthm\label{bdd1}
Let $\psi$ be a complex-valued function on $T$  and let $0<p<\infty$. Then the following are equivalent:
\begin{enumerate}
\item[(a)] The multiplication operator $M_{\psi}:\mathbb{T}_p\rightarrow \mathbb{T}_\infty$ is bounded.

\item[(b)] The multiplication operator $M_{\psi}:\mathbb{T}_{p,0}\rightarrow \mathbb{T}_{\infty,0}$ is bounded.

\item[(c)] $\ds\sup\limits_{n\in \mathbb{N}_0}\left\{c_n^\frac{1}{p}\, M_\infty(n,\psi) \right\}<\infty$.
\end{enumerate}
Moreover,  $\|M_\psi\|=\sup\limits_{n\in \mathbb{N}_0}\left\{c_n^\frac{1}{p} M_\infty(n,\psi) \right\}$.
\ethm

\bpf
For $f\in \mathbb{T}_p$ and $n\in \N_0$, we have
\beqq
M_\infty(n,M_\psi f) &=& \sup_{|v|=n}|\psi(v)|\,|f(v)|\\
                   &\leq & M_\infty(n,\psi)\, M_\infty(n,f)\\
                   &\leq & c_n^\frac{1}{p}M_\infty(n,\psi) M_p(n,f) \mbox {~\,\,\,(by Lemma \ref{inequ}).~}
\eeqq
This gives that $(c)\Rightarrow (a)$ and  $(c)\Rightarrow (b)$ with
$$
\|M_\psi\|\leq \sup\limits_{n\in \mathbb{N}_0}\left\{c_n^\frac{1}{p} M_\infty(n,\psi) \right\}.
$$
For the other way implications, fix $n\in \mathbb{N}_0$ and choose $v_n\in D_n$ such that 
$$\max\limits_{|v|=n}|\psi(v)|=|\psi(v_n)|.$$
Take $f=(c_n)^\frac{1}{p}\chi_{v_n}$ so that $\|f\|_p=1$ and

$$M_\infty(n,\psi f)=\max_{|v|=n}|\psi(v)|\,|f(v)|=|\psi(v_n)|\,|f(v_n)|=c_n^\frac{1}{p} M_\infty(n,\psi).$$
Thus,
$$c_n^\frac{1}{p} M_\infty (n,\psi)= M_\infty(n,\psi f)= \|\psi f \|_\infty \leq \|M_\psi \|. $$
Since $n\in \mathbb{N}_0$ was arbitrary and $f$ is in $\mathbb{T}_{p,0}$ and as well in $\mathbb{T}_{p}$, we get the implications $(b)\Rightarrow (c)$ and  $(a)\Rightarrow (c)$ with
$$
\|M_\psi\|\geq \sup\limits_{n\in \mathbb{N}_0}\left\{c_n^\frac{1}{p} \, M_\infty(n,\psi) \right\}.
$$
The theorem follows.
\epf

\bthm
Let $\psi$ be a function from $T$ to $\mathbb{C}$ and let $0<p<q<\infty$. Then the following are equivalent:
\begin{enumerate}
\item[(a)] The multiplication operator $M_{\psi}:\mathbb{T}_p\rightarrow \mathbb{T}_q$ is bounded.

\item[(b)] The multiplication operator $M_{\psi}:\mathbb{T}_{p,0}\rightarrow \mathbb{T}_{q,0}$ is bounded.

\item[(c)] $\ds\sup\limits_{n\in \mathbb{N}_0}\left\{c_n^{\frac{1}{p}-\frac{1}{q}} M_\infty(n,\psi)   \right\} < \infty$.
\end{enumerate}
In this case, $$\|M_\psi\|= \sup_{n\in \mathbb{N}_0} \left\{c_n^{\frac{1}{p}-\frac{1}{q}} M_\infty(n,\psi)   \right\}. $$
\ethm

\bpf
For $f\in \mathbb{T}_p$ and $n\in \N_0$, we have
$$M_q(n,\psi f)\leq M_\infty(n,\psi)\, M_q(n,f)\leq c_n^{\frac{1}{p}-\frac{1}{q}}M_\infty(n,\psi)\, M_p(n,f) \mbox {~\,\,\,(by Lemma \ref{inequ}).~}$$
That is,
$$M_q(n,\psi f)\leq \left\{c_n^{\frac{1}{p}-\frac{1}{q}}M_\infty(n,\psi) \right\} M_p(n,f) \mbox{~for each $n\in \N_0$}.$$
This yields that $(c)\Rightarrow (a)$ and  $(c)\Rightarrow (b)$ with
$$
\|M_\psi\|\leq \sup\limits_{n\in \mathbb{N}_0}\left\{c_n^{\frac{1}{p}-\frac{1}{q}} M_\infty(n,\psi)\right\}.
$$
For the converse implications, fix $n\in \mathbb{N}_0$ and
choose $v_n\in D_n$ such that $M_\infty(n,\psi)=|\psi(v_n)|$.
Take $f=\chi_{\{v_n\}} c_n^\frac{1}{p}$ so that  $\|f\|_p=1$ and
$$\|\psi f \|_q =M_q(n,\psi f)=\left(\frac{1}{c_n} |\psi(v_n)|^q |f(v_n)|^q \right)^\frac{1}{q}=c_n^{-\frac{1}{q}}|\psi(v_n)|\,|f(v_n)|=c_n^{\frac{1}{p}-\frac{1}{q}} M_\infty(n,\psi) $$
and thus,
$$\|M_\psi\|=\sup_{\|f\|_p =1}\left\{\|\psi f \|_q \right\}\geq \sup_{n\in \mathbb{N}_0} \left\{c_n^{\frac{1}{p}-\frac{1}{q}} \, M_\infty(n,\psi)   \right\}.$$
Since the functions $f$  used in this part of proof lies in $\mathbb{T}_{p,0}$ as well in $\mathbb{T}_{p}$, it proves the implications $(b)\Rightarrow (c)$ and  $(a)\Rightarrow (c)$.
The norm estimate
$$\|M_\psi\|= \ds\sup_{n\in \mathbb{N}_0} \left\{c_n^{\frac{1}{p}-\frac{1}{q}}\, M_\infty(n,\psi)   \right\} $$
follows from the lines of the proof.
\epf

\section{Isometry}\label{Msha1-Sec4}

A linear operator $A$ from a Banach space $(X,\|.\|_X)$ to a Banach space $(Y,\|.\|_Y)$ is said to be an isometry if $\|Ax\|_Y = \|x\|_X$ for all $x\in X$.

In this section, we discuss about isometric multiplication operators between various $\mathbb{T}_p$ spaces. First, we recall a result from \cite{MP-Tp-spaces} that characterizes isometric multiplication operators on $\mathbb{T}_p$.

\begin{Thm} \rm{(\cite[Theorem 4.10]{MP-Tp-spaces})}
Let $X$ be either $\mathbb{T}_{p}$ or $\mathbb{T}_{p,0}$, where $0<p\leq\infty$ and let $M_\psi: X \rightarrow X$ be a bounded operator on $X$. Then
$M_\psi$ is an isometry on $X$ if and only if $|\psi(v)|=1$ for all $v\in T$.
\end{Thm}
\bthm\label{iso}
Let $T$ be a rooted tree such that $\{c_n\}$ is unbounded and let $p,q\in (0,\infty]$ with $p\neq q$. Then, there is no isometric multiplication operator $M_\psi$ from $\mathbb{T}_p$ to $\mathbb{T}_q$.
\ethm
\bpf
We prove our claim case by case.\\
\\
\textbf{Case 1:} ($M_\psi:\mathbb{T}_p\rightarrow \mathbb{T}_\infty$, where $0<p< \infty$)

Fix $v\in T$ and take $f=(c_{|v|})^\frac{1}{p}\chi_{v}$, so that $\|f\|_p=1$.
For $M_\psi$ to be an isometry, we should have $\|\psi f \|_\infty=1$, which implies $(c_{|v|})^\frac{1}{p}\,|\psi(v)|=1$. Since $v\in T$ was arbitrary, we get $|\psi(v)|=(c_{|v|})^{-\frac{1}{p}}$ for all $v\in T$.

Choose $k\in \N$ such that $c_k>1$ and define
$$ f(v)=\left\{
\begin{array}{cl}
1& \mbox{ if } v\in D_k\\
0 & \mbox{ otherwise}.
\end{array}
\right.
$$
Then $\|f\|_p=1$, but
$$\|\psi f \|_\infty= \sup_{v\in T}|\psi(v)\|f(v)|=\sup_{v\in D_k}|\psi(v)|=c_k^{-\frac{1}{p}}\neq 1.$$
This forces that $M_\psi: \mathbb{T}_p\rightarrow \mathbb{T}_\infty$ cannot be an isometry.\\
\\
\textbf{Case 2:} ($M_\psi:\mathbb{T}_\infty \rightarrow \mathbb{T}_p$, where $0<p< \infty$)

For each $v\in T$, we have $\|\chi_{v}\|_\infty=1$.
For $M_\psi$ to be an isometry, we must have
$$\|\psi \chi_{v} \|_p=(c_{|v|})^{-\frac{1}{p}}|\psi(v)|=1.$$
This yields that
$|\psi(v)|=(c_{|v|})^\frac{1}{p}$ for all $v\in T$. Therefore, for each $n\in \N_0$, we have $M_p(n,\psi)=c_n^\frac{1}{p}$, which implies that $\psi\notin \mathbb{T}_p$. Then by Theorem \ref{bdd2}, $M_\psi$ is not a bounded operator. Hence, there does not exist an isometric multiplication operator from $\mathbb{T}_\infty$ to $\mathbb{T}_p$.\\
\\
\textbf{Case 3:} ($M_\psi:\mathbb{T}_p\rightarrow \mathbb{T}_q$ with $0<q<p<\infty$)

For each $v\in T$, define $f=(c_{|v|})^\frac{1}{p}\,\chi_{v}$ so that $\|f\|_p=1$.
For $M_\psi$ to be an isometry, we should have $\|\psi f \|_q=1$ which implies that
$$|\psi(v)|=(c_{|v|})^{\frac{1}{q}-\frac{1}{p}} \mbox{~for all~} v\in T.$$
Since $\frac{1}{q}-\frac{1}{p}>0$ and $\{c_n\}$ is unbounded, the sequence  $\left\{M_{\frac{pq}{p-q}}(n,\psi) \right\}$, i.e., $\left\{c_n^{\frac{1}{q}-\frac{1}{p}}\right\}$, is not bounded.
Thus, $M_\psi$ is a unbounded operator by Theorem \ref{thm1}. Hence,
$M_\psi:\mathbb{T}_p\rightarrow\mathbb{T}_q$ cannot be isometry.\\
\\
\textbf{Case 4:} ($M_\psi:\mathbb{T}_p\rightarrow \mathbb{T}_q$ with $0<p<q<\infty$)

By similar arguments as in Case 3, we have
$$|\psi(v)|=(c_{|v|})^{\frac{1}{q}-\frac{1}{p}} \mbox{~for all~} v\in T.$$ Choose $v_1,v_2\in T$ such that $|v_1|=|v_2|\,(= k),$ say. Take $f=\chi_{\{v_1,v_2\}}$, that is,
$$ f(v)=\left\{
\begin{array}{cl}
1& \mbox{~if~} v=v_1 \mbox{~or~} v=v_2\\
0 & \mbox{ otherwise}.
\end{array}
\right.
$$
Then, we have
$$\|f\|_p=\left(\frac{2}{c_k} \right)^\frac{1}{p} \mbox{~and~}
\|\psi f \|_q=\left(\frac{2 c_k^{(1-\frac{q}{p})}}{c_k}  \right)^\frac{1}{q}=\frac{2^\frac{1}{q}}{c_k^\frac{1}{p}}. $$

If $p\neq q$, then $2^\frac{1}{p}\neq 2^\frac{1}{q}$ and thus $\|f\|_p\neq \|\psi f \|_q$. Hence $M_\psi:\mathbb{T}_p\rightarrow \mathbb{T}_q$ fails to be an isometry.
\epf

Test functions that we used in Theorem \ref{iso}, lie in $\mathbb{T}_{p,0}$. Therefore, we have the following result.
\bcor
Let $T$ be a rooted tree such that $\{c_n\}$ is unbounded and let $p,q\in (0,\infty]$ with $p\neq q$. Then, there does not exist an  isometric multiplication operator $M_\psi$ from $\mathbb{T}_{p,0}$ to $\mathbb{T}_{q,0}$.
\ecor
\section{Compact Multiplication Operators}\label{Msha1-Sec5}
This section is devoted to compact multiplication operators between various  $\mathbb{T}_p$ spaces. In view of \cite[Theorem 7]{CO-Tp-spaces} and \cite[Lemma 4.8]{MP-Tp-spaces}, we prove a general result. The following result characterizes compact operators between $\mathbb{T}_p$ spaces.
\bthm
Let $T$ be a rooted tree and let $A$ be a linear operator such that $Af_n\rightarrow 0$ pointwise in  $\mathbb{T}_q$ whenever  $f_n\rightarrow 0$  pointwise in  $\mathbb{T}_p$.
Then, $A:\mathbb{T}_p\rightarrow \mathbb{T}_q$ is compact  if and only if
for every bounded sequence $\{f_n\}$ in  $\mathbb{T}_p$ that converges to $0$ pointwise, the sequence $\|A f_n \|_q\rightarrow 0 $ as $n\rightarrow \infty$.

\ethm

\bpf
Assume that $A:\mathbb{T}_p\rightarrow \mathbb{T}_q$ is compact.  Suppose that $\{f_n\}$  is a bounded sequence in $\mathbb{T}_p$, converging to 0 pointwise. Then, there is a subsequence $\{f_{n_k}\}$ of $\{f_n\}$ such that $\{A f_{n_k}\}$ converges in $\|.\|_q$ to some function, say, $f$. It follows that $\{A f_{n_k}\}$ converges to $f$ pointwise. Since the pointwise convergence of $\{f_n\}$ to 0 implies that $f\equiv 0$. Thus, we deduce that $\{A f_{n_k}\}$ converges to 0 in $\|.\|_q$.
Now, it is easy to verify that the sequence $\|A f_n \|_q\rightarrow 0 $ as $n\rightarrow \infty$.


To prove the converse part, let us begin by assuming that for every bounded sequence $\{f_n\}$ in  $\mathbb{T}_p$ converging to 0 pointwise, then the sequence $\|A f_n \|_q\rightarrow 0 $ as $n\rightarrow \infty$. Let $\{g_n\}$ be a bounded sequence in  $\mathbb{T}_p$. Without loss of generality, assume that $\|g_n\|_p\leq 1$ for all $n$. Then, by Lemma \Ref{lem:bound}, for each fixed $v\in T$, $\{g_n(v)\}$ is a bounded sequence. By the diagonalization process, there is a subsequence $\{g_{n_k}\}$ of $\{g_n\}$ such that $\{g_{n_k}\}$ converges pointwise to $g$ (say). Since $M_p(m,g_{n_k})\rightarrow M_p(m,g)$ for each $m\in \mathbb{N}_0$ and $\|g_{n_k}\|\leq 1$ for all $k$, we have $M_p(m,g)\leq 1$ for all $m$ which implies that $g\in\mathbb{T}_p$ with $\|g\|_p\leq 1$.
Take $f_k=g_{n_k}-g\in \mathbb{T}_p$ so that $f_k$ converges to $0$ pointwise. By our assumption, the sequence $\|A f_k \|_q\rightarrow 0 $ as $k\rightarrow \infty$. This proves that
$A$ is a compact operator.
\epf

Since $\psi f_n\rightarrow 0$ pointwise whenever $f_n\rightarrow 0$ pointwise for any complex-valued map $\psi$, we have a following consequence.

\bcor\label{compactthm}
If $p,q\in (0,\infty]$, then the multiplication operator $M_\psi:\mathbb{T}_p\rightarrow \mathbb{T}_q$ is compact  if and only if every bounded sequence $\{f_n\}$ in  $\mathbb{T}_p$ converging to $0$ pointwise, the sequence $\|\psi f_n \|_q\rightarrow 0 $ as $n\rightarrow \infty$.
\ecor

\bcor
Let $p,q\in (0,\infty]$ and the multiplication operator $M_\psi:\mathbb{T}_{p,0}\rightarrow \mathbb{T}_{q,0}$ be compact. Then, for every bounded sequence $\{f_n\}$ in  $\mathbb{T}_{p,0}$ converging to $0$ pointwise, the sequence $\|\psi f_n \|_q\rightarrow 0 $ as $n\rightarrow \infty$.
\ecor
The following simple observations will be useful in the sequal.
\begin{enumerate}
\item For any complex sequence $\{x_n\}$, $x_n \rightarrow 0$ as $n\rightarrow \infty$ if and only if $\sup\limits_{n\geq N}|x_n|\rightarrow 0$ as $N\rightarrow \infty$.

\item Finite rank operators are compact and norm limit of a sequence of compact operators is again a compact operator.
\end{enumerate}

Next, we discuss the compactness of multiplication operator between $\mathbb{T}_{p}$ spaces case by case.

\begin{Thm}\rm{\cite[Theorem 4.7]{MP-Tp-spaces}}
Let $X$ be either $\mathbb{T}_{p}$ or $\mathbb{T}_{p,0}$, where $0<p\leq\infty$ and let $M_\psi: X \rightarrow X$ be a bounded multiplication operator on $X$. Then
$M_\psi$ is a compact operator on $X$ if and only if $\psi(v)\rightarrow 0$ as $|v|\rightarrow \infty$.
\end{Thm}
\bthm\label{cpt}
Let $M_\psi:\mathbb{T}_p\rightarrow \mathbb{T}_\infty$ $(0<p<\infty)$ be a bounded multiplication operator. Then $M_\psi$ is a compact operator if and only if
$$
c_n^{\frac{1}{p}} M_\infty (n,\psi)\rightarrow 0 \mbox{~as~ $n\rightarrow \infty$, i.e.,~} M_\infty(n,\psi)=o(c_n^{-\frac{1}{p}}).
$$
\ethm
\bpf
Let $\psi$ be an arbitrary function such that $c_n^{\frac{1}{p}} M_\infty (n,\psi)\rightarrow 0$ as $n\rightarrow \infty$. Consider $\{\psi_n\}_{n\in \N}$, where
$$ \psi_n(v)=\left\{
\begin{array}{cl}
\psi(v)& \mbox{ if } |v|\leq n\\
0 & \mbox{ otherwise}.
\end{array}
\right.
$$
Thus, $M_{\psi_n}$ is a finite rank operator and hence, is a compact operator for every $n$. Moreover,
$$\|M_{\psi_n} - M_\psi \|=\|M_{\psi_n-\psi}\|=\sup_{k\in \mathbb{N}_0}\left\{c_k^\frac{1}{p} M_\infty
(k,\psi_n -\psi) \right\}=\sup_{k> n}\left\{c_k^\frac{1}{p} M_\infty
(k,\psi) \right\}, $$
which approaches to zero as $n\rightarrow \infty$. Hence, $M_\psi$, being the limit of the  sequence $\{M_{\psi_n}\}$ of compact operators, is a compact operator.

Conversely, assume that $M_\psi$ is compact. For each $n\in\N_0$, choose $v_n\in D_n$ such that $|\psi(v_n)|=M_\infty (n,\psi)$ and define $f_n$ by
$$ f_n(v)=\left\{
\begin{array}{cl}
c_n^\frac{1}{p}\chi_{v_n}& \mbox{ if } v=v_n\\
0 & \mbox{ otherwise}.
\end{array}
\right.
$$
Then, $\|f_n\|_p=1$ for all $n\in \mathbb{N}_0$ and for each fixed $v\in T$,  $f_n(v)=0$ for all $n>|v|$. This implies that $f_n$ converges pointwise to $0$. Therefore, by Corollary \ref{compactthm}, we get
$$\|\psi f_n \|_\infty=c_n^{\frac{1}{p}}|\psi(v_n)|=c_n^{\frac{1}{p}}M_\infty (n,\psi)\rightarrow 0 ~~\text{as}~~ n\rightarrow \infty. $$
The desired result follows.
\epf
\bcor
For $0<p<\infty$, $M_\psi:\mathbb{T}_{p,0}\rightarrow \mathbb{T}_{\infty,0}$  is a compact operator if and only if
$$
c_n^{\frac{1}{p}} M_\infty (n,\psi)\rightarrow 0 \mbox{~as~ $n\rightarrow \infty$, i.e.,~} M_\infty(n,\psi)=o(c_n^{-\frac{1}{p}}).
$$
\ecor

\bthm
Let $M_\psi:\mathbb{T}_\infty\rightarrow \mathbb{T}_p$ $(0<p<\infty)$ be a bounded multiplication operator. Then $M_\psi$ is a compact operator if and only if $M_p(n,\psi)\rightarrow 0$ as $n\rightarrow \infty$, i.e., $\psi\in \mathbb{T}_{p,0}$.
\ethm

\bpf
Suppose that $ M_p (n,\psi)\rightarrow 0$ as $n\rightarrow \infty$. Define $\{\psi_n\}$ as in Theorem \ref{cpt}.
Then, $M_{\psi_n}$ is a compact operator for every $n$. Since
$$\|M_{\psi_n} - M_\psi \|=\|\psi_n-\psi\|_p=\sup_{k\in \mathbb{N}_0} M_p(k,\psi_n -\psi)=\sup_{k> n} M_p
(k,\psi) $$
approaches to zero as $n\rightarrow \infty$,  $M_\psi$ is a compact operator.

Convesely, assume that $M_\psi$ is a compact  operator. Define $f_n$ by
$$ f_n(v)=\left\{
\begin{array}{cl}
1 & \mbox{ if } v\in D_n\\
0 & \mbox{ otherwise}.
\end{array}
\right.
$$
Then, $\|f_n\|_\infty=1$ for all $n\in \mathbb{N}_0$ and $f_n$ converges to $0$ pointwise. Corollary \ref{compactthm} gives that,
$$\|\psi f_n \|_p=M_p (n,\psi)\rightarrow 0 ~~\text{as}~~ n\rightarrow \infty. $$
It completes the proof.
\epf

\bcor For $0<p<\infty$,
 $M_\psi:\mathbb{T}_{\infty,0}\rightarrow \mathbb{T}_{p,0}$ is a compact operator if and only if $M_p(n,\psi)\rightarrow 0$ as $n\rightarrow \infty$, i.e., $\psi\in \mathbb{T}_{p,0}$.
\ecor

\bthm
Let $M_\psi:\mathbb{T}_p\rightarrow \mathbb{T}_q$ $(0<q<p<\infty)$ be a bounded multiplication operator. Then $M_\psi$ is a compact operator if and only if
$$ M_{\frac{pq}{p-q}}(n,\psi)\rightarrow 0 \mbox{~as~ $n\rightarrow \infty$, i.e.,~} \psi \in \mathbb{T}_{\frac{pq}{p-q},0}.$$
\ethm

\bpf
Let $\psi$ be an arbitrary function such that $M_{\frac{pq}{p-q}}(n,\psi)\rightarrow 0$ as $n\rightarrow \infty$. For each $n$, define $\{\psi_n\}$ as in Theorem \ref{cpt}.
Then, $M_{\psi_n}$ is a compact operator for every $n$ and
$$\|M_{\psi_n}-M_\psi \|=\|\psi_n -\psi\|_{\frac{pq}{p-q}}=\sup_{k\in \mathbb{N}}M_{\frac{pq}{p-q}}(k,\psi_n -\psi) =\sup_{k> n} M_{\frac{pq}{p-q}} (k,\psi).$$
By the hypothesis, $\|M_{\psi_n}-M_\psi \|\rightarrow0$ as $k\rightarrow \infty$. Hence, $M_\psi$ is a compact operator.

For the proof of converse part, let us assume that $M_\psi$ is compact. Define $f_n$ by
$$ f_n(v)=\left\{
\begin{array}{cl}
\frac{1}{A_n} |\psi(v)|^{\frac{q}{p-q}}& \mbox{ if } v\in D_n\\
0 & \mbox{ otherwise},
\end{array}
\right.
$$
where $A_n=\ds\left(\frac{1}{c_n} \sum\limits_{|v|=n}|\psi(v)|^{\frac{pq}{p-q}} \right)^\frac{1}{p}$ so that $\|f_n\|_p=1$ for all $n\in \mathbb{N}_0$. It is easy to see that $f_n$ converges to $0$ pointwise and

$$
\begin{array}{rcl}
\|\psi f_n \|_q& = & M_q(n,\psi f_n)\\
                &= &\ds \left\{\frac{1}{c_n}\sum\limits_{|v|=n}\left(|\psi(v)|\frac{|\psi(v)|}{A_n}^{\frac{q}{p-   q}} \right)^q \right\}^\frac{1}{q} \\
                &= &\ds \frac{1}{A_n}\left(\frac{1}{c_n}\sum\limits_{|v|=n}|\psi(v)|^{\frac{pq}{p-q}} \right)^\frac{1}{q} \\
                &=& M_{\frac{pq}{p-q}}(n,\psi).
\end{array}
$$
By Corollary \ref{compactthm}, we have
$M_{\frac{pq}{p-q}}(n,\psi) \rightarrow 0 ~~\text{as}~~ n\rightarrow \infty.$ Hence the theorem.
\epf

\bcor
Let $M_\psi:\mathbb{T}_{p,0}\rightarrow \mathbb{T}_{q,0}$ $(0<q<p<\infty)$ be a bounded multiplication operator. Then $M_\psi$ is a compact operator if and only if
$$ M_{\frac{pq}{p-q}}(n,\psi)\rightarrow 0 \mbox{~as~ $n\rightarrow \infty$, i.e.,~} \psi \in \mathbb{T}_{\frac{pq}{p-q},0}.$$
\ecor

\bthm
Let $M_\psi:\mathbb{T}_p\rightarrow \mathbb{T}_q$ $(0<p<q<\infty)$ be a bounded multiplication operator. Then $M_\psi$ is a compact operator if and only if
$$c_n^{\frac{1}{p}-\frac{1}{q}} M_\infty (n,\psi)\rightarrow 0\mbox{~as~ $n\rightarrow \infty$,  i.e.,~} M_\infty(n,\psi)=o(c_n^{\frac{1}{q}-\frac{1}{p}}).
$$
\ethm
\bpf
Let $\psi$ be an arbitrary complex-valued function such that $c_n^{\frac{1}{p}-\frac{1}{q}} M_\infty (n,\psi)\rightarrow 0$ as $n\rightarrow \infty$. For each $n$, define $\{\psi_n\}$ as in Theorem \ref{cpt},
so that $M_{\psi_n}$ is a compact operator for every $n$. Moreover,
$$\|M_{\psi_n} - M_\psi \|=\|M_{\psi_n-\psi}\|=\sup_{k\in \mathbb{N}_0}\left\{c_k^{\frac{1}{p}-\frac{1}{q}} M_\infty
(k,\psi_n -\psi) \right\}=\sup_{k> n}\left\{c_k^{\frac{1}{p}-\frac{1}{q}} M_\infty
(k,\psi) \right\}, $$
which converges to zero as $n\rightarrow \infty$. This forces that, $M_\psi$ is compact.

For the proof of converse part, assume that  $M_\psi$ is a compact operator. For $n\in\N_0$, choose $v_n\in D_n$ in such a way that $|\psi(v_n)|=M_\infty (n,\psi)$ and define $f_n$ by
$$ f_n(v)=\left\{
\begin{array}{cl}
c_n^\frac{1}{p}\chi_{v_n}& \mbox{ if } v=v_n\\
0 & \mbox{ otherwise}.
\end{array}
\right.
$$
Thus, we have $\|f_n\|_p=1$ for all $n\in \mathbb{N}_0$ and $f_n$ converges pointwise to $0$. As a consequence of Corollary \ref{compactthm}, we see that
$$\|\psi f_n \|_q={c_n^{-\frac{1}{q}}}|f_n(v_n)\|\psi(v_n)|=c_n^{\frac{1}{p}-\frac{1}{q}}M_\infty(n,\psi)\rightarrow 0 ~~\text{as}~~ n\rightarrow \infty. $$
This proves our claim.
\epf
\bcor
Let $M_\psi:\mathbb{T}_{p,0}\rightarrow \mathbb{T}_{q,0}$ $(0<p<q<\infty)$ be a bounded multiplication operator. Then $M_\psi$ is a compact operator if and only if
$$c_n^{\frac{1}{p}-\frac{1}{q}} M_\infty (n,\psi)\rightarrow 0\mbox{~as~ $n\rightarrow \infty$ i.e.,~} M_\infty(n,\psi)=o(c_n^{\frac{1}{q}-\frac{1}{p}}).
$$
\ecor

\section{Invertible Multiplication Operators}\label{Msha1-Sec6}


We devote this section for the discussion on invertiblity of multiplication operators between various $\mathbb{T}_p$ spaces,

\bprop\label{inv}
Let $p,q\in (0,\infty]$ and $\psi:T\rightarrow \mathbb{C}$ be an arbitrary function. If $M_{\psi}:\mathbb{T}_p\rightarrow \mathbb{T}_q$ is invertible, then $\psi(v)\neq 0$  for all $v\in T$ and $(M_\psi)^{-1}=M_{\frac{1}{\psi}}:\mathbb{T}_q\rightarrow \mathbb{T}_p$.
\eprop

\bpf
Suppose $M_\psi$ is an invertible operator from $\mathbb{T}_p$ to $\mathbb{T}_q$. Then, there exists a bounded linear operator $B:\mathbb{T}_q\rightarrow \mathbb{T}_p$ such that
$$M_\psi(Bf)=f \mbox{~~for all~~} f\in \mathbb{T}_q.$$
In particular, $\psi(v)Bf(v)=f(v)$ for all $v\in T$ and for all $f\in \mathbb{T}_q$. For each fixed $v\in T$, by setting $f=\chi_v$, characteristic function on $\{v\}$, we see that $\psi(v)Bf(v)=1$ and thus $\psi(v)\neq 0$. Since $v\in T$ was arbitrary, $\psi(v)\neq 0$ for all $v\in T$ and thus,
$$Bf(v)=\frac{1}{\psi(v)} f(v) \mbox{~ for all~} v\in T.$$
Therefore, $B=A^{-1}$ is a multiplication operator induced by $\frac{1}{\psi}$. It completes the proof.
\epf

\bthm
Let $M_\psi:\mathbb{T}_p\rightarrow \mathbb{T}_p$ $(0<p\leq \infty)$ be a bounded multiplication operator on $\mathbb{T}_p$. Then $M_\psi$ is invertible on $\mathbb{T}_p$ if and only if there exist $m,M>0$ such that $m\leq |\psi(v)|\leq M$ for all $v\in T$.
\ethm

\bpf
Suppose that $M_\psi$ is invertible on $\mathbb{T}_p$. Then, by Theorem \Ref{bdd}, $\psi$ and $\frac{1}{\psi}$ are bounded functions. Hence, there exist $m,M>0$ such that $m\leq |\psi(v)|\leq M$ for all $v\in T$.

Conversely, assume that there exist $m,M>0$ such that $m\leq |\psi(v)|\leq M$ for all $v\in T$. We have $\left|\frac{1}{\psi(v)} \right|\leq \frac{1}{m}<\infty$ for all $v\in T$. Therefore $M_\frac{1}{\psi}$
is bounded operator on $\mathbb{T}_p$ by Theorem \Ref{bdd}. Moreover,
$$M_{\frac{1}{\psi}} (M_{\psi}(f))=M_\psi(M_{\frac{1}{\psi}} (f))=f \mbox{~~for all~~} f\in \mathbb{T}_p.$$
Thus, $M_\psi$ is an invertible operator on $\mathbb{T}_p$. The desired result follows.
\epf
\bcor
Let $T$ be a rooted tree and $0<p\leq \infty$. Then, $M_\psi:\mathbb{T}_{p,0}\rightarrow \mathbb{T}_{p,0}$  is invertible  if and only if there exist $m,M>0$ such that $m\leq |\psi(v)|\leq M$ for all $v\in T$.
\ecor

\bthm\label{inv2}
Let $0<p<\infty$ and $T$ be a rooted tree such that $\{c_n\}$ is unbounded. Then, there are no invertible multiplication operators from $\mathbb{T}_\infty$ to $\mathbb{T}_p$.
\ethm
\bpf
Suppose there exists  $\psi:T\rightarrow \mathbb{C}$ such that $M_\psi:\mathbb{T}_\infty\rightarrow \mathbb{T}_p$ is an invertible operator. Then, by Proposition \ref{inv},
$\psi(v)\neq 0$ for every $v\in T$ and $M_{\frac{1}{\psi}}:\mathbb{T}_p\rightarrow \mathbb{T}_\infty$ is a bounded operator. Theorem \ref{bdd1} assures the existance of a positive constant $K>0$ such that
$$
M_\infty(n,\frac{1}{\psi})c_n^{\frac{1}{p}} \leq K \mbox{~for all~} n\in \mathbb{N}_0.
$$
Upon setting $m_\infty(n,\psi)=\min\limits_{|v|=n}|\psi(v)|$, we see that
$$\frac{c_n^{\frac{1}{p}}}{m_\infty(n,\psi)}\leq K \mbox{~for all~} n\in \mathbb{N}_0.
$$
Fix $n\in \mathbb{N}_0$. Since $m_\infty(n,\psi)\leq |\psi(v)|\mbox{~for all~} v\in D_n$, we have
$$\frac{c_n^{\frac{1}{p}}}{K}\leq |\psi(v)|\mbox{~for all~} v\in D_n.$$
This implies that
$$\frac{c_n^{\frac{1}{p}}}{K}\leq M_p(n,\psi) \mbox{~for all~} n\in \mathbb{N}_0. $$
Since $\{c_n\}$ is unbounded, the sequence $\{M_p(n,\psi)\}$ has to be unbounded and thus, $\psi\notin \mathbb{T}_p$ which in turn implies that $\psi$ cannot induce a bounded multiplication operator from $\mathbb{T}_\infty$ to $\mathbb{T}_p$ by Theorem \ref{bdd2}.  Hence, there are no invertible multiplication operators from $\mathbb{T}_\infty$ to $\mathbb{T}_p$ as desired.
\epf

If $M_\psi:\mathbb{T}_p\rightarrow \mathbb{T}_\infty$ is an invertible operator then $M_{\frac{1}{\psi}}:\mathbb{T}_\infty \rightarrow \mathbb{T}_p$ is invertible. This observation gives the following corollary.

\bcor
Let $0<p<\infty$ and $T$ be a rooted tree such that $\{c_n\}$ is unbounded. Then, there exists no invertible multiplication operator $M_\psi:\mathbb{T}_p \rightarrow \mathbb{T}_\infty$.
\ecor

\bthm
Let $p,q\in (0,\infty)$ with $p\neq q$ and $T$ be a rooted tree such that $\{c_n\}$ is unbounded. Then, there is no invertible multiplication operators from $\mathbb{T}_p$ to $\mathbb{T}_q$.
\ethm
\bpf
We omit the detail, since proof is similar to that of Theorem \ref{inv2}.
\epf

\bcor
Let $p,q\in (0,\infty]$ with $p\neq q$ and $T$ be a rooted tree such that $\{c_n\}$ is unbounded. Then, there does not exist an invertible multiplication operator $M_\psi:\mathbb{T}_{p,0}\rightarrow \mathbb{T}_{q,0}$.
\ecor

\bthm\label{inj}
Let $p,q\in(0,\infty]$ and let $\psi:T\rightarrow \mathbb{C}$ be an arbitrary function. Then the multiplication operator $M_\psi:\mathbb{T}_p\rightarrow \mathbb{T}_q$ is injective if and only if $\psi(v)\neq 0$ for all $v\in T$.
\ethm
\bpf
Suppose that $\psi(v)\neq 0$ for all $v\in T$ and $M_\psi f=M_\psi g$ for some $f,g\in \mathbb{T}_p$. This leads that $\psi(v)(f-g)(v)=0$ for all $v\in T$. Therefore, $f(v)=g(v)$ for all $v\in T$ and hence, $M_\psi$ is injective.

Conversely, assume that $M_\psi$ is injective. Suppose that $\psi(v_0)=0$ for some $v_0\in T$. Take $f=\chi_{v_0}$ so that $M_\psi f(v)=\psi(v)\chi_{v_0}(v)=0$ for all $v\in T$. Thus, there exists a function $f\neq 0$ such that $M_\psi f=0$, which contradicts the fact that, $M_\psi$ is injective. Hence $\psi(v)\neq 0$ for all $v\in T$.
\epf
\bprop\label{bij}
If $M_\psi:\mathbb{T}_p\rightarrow \mathbb{T}_q$ is surjective, then  $M_\psi:\mathbb{T}_p\rightarrow \mathbb{T}_q$ is injective.
\eprop
\bpf
Suppose that $M_\psi:\mathbb{T}_p\rightarrow \mathbb{T}_q$ is surjective. Fix $w\in T$. Then, as    $\chi_w \in \mathbb{T}_q$, there exist $f\in \mathbb{T}_p$ such that $M_\psi(f)=\chi_w$. This means that $\psi(v)f(v)=\chi_w(v)$ for all $v\in T$. In particular $\psi(w)f(w)=1$ and therefore, $\psi(w)\neq 0$. Since $w\in T$ was arbitrary, we get $\psi(v)\neq 0$ for all $v\in T$.  Hence, by Theorem \ref{inj}, $M_\psi:\mathbb{T}_p\rightarrow \mathbb{T}_q$ is injective.
\epf
\bthm
Let $T$ be a rooted tree such that $\{c_n\}$ is unbounded and let $p,q\in [1,\infty]$ with $p\neq q$. Then no bounded multiplication operator $M_\psi:\mathbb{T}_p\rightarrow \mathbb{T}_q$ is onto.
\ethm
\bpf
Suppose that $M_\psi:\mathbb{T}_p\rightarrow \mathbb{T}_q$ is onto.  Then, by Proposition \ref{bij}, $M_\psi:\mathbb{T}_p\rightarrow \mathbb{T}_q$ is bijective. Note that $\mathbb{T}_r$ is a Banach space for each $r\geq1$. By inverse mapping theorem \cite[Theorem 12.5]{Conway:Book}, $M_\psi:\mathbb{T}_p\rightarrow \mathbb{T}_q$ is an invertible operator, which is a contradiction to the fact that there does not exist an invertible multiplication operator between $\mathbb{T}_p$ and $\mathbb{T}_q$ when $p\neq q$. Thus, a bounded multiplication operator $M_\psi:\mathbb{T}_p\rightarrow \mathbb{T}_q$ never be onto.
\epf

\section{Fixed points}\label{Msha1-Sec7}
Let $A$ be a linear operator on a normed linear space $X$. An element $x\in X$ is called a fixed point of $A$ if $Ax=x$.

\bthm
Let $M_\psi:\mathbb{T}_p \rightarrow \mathbb{T}_p$ be a bounded multiplication operator and let $E=\{v\in T: \psi(v)\neq 1 \}$. Then, the collection of all fixed points of $M_\psi$ is preciously the closed subspace $\{f\in \mathbb{T}_p:f|_{E}=0\}$.
\ethm

\bpf
Since $M_\psi$ is bounded operator on $\mathbb{T}_p$, we have $\psi\in \mathbb{T}_\infty$ and thus $\psi-1\in \mathbb{T}_\infty$, which shows that $M_{\psi-1}$ is a bounded operator on $\mathbb{T}_p$. A function $f\in \mathbb{T}_p$ is a fixed point of $M_\psi$ if and only if $(\psi-1)f=0$, i.e., $f\in \text{ker}(M_{\psi -1})$. Since kernel of a bounded operator is a closed subspace, the collection of all fixed points of $M_\psi$ is a closed subspace of $\mathbb{T}_p$ so that the conclusion
$(\psi-1)f=0$ forces that $f=0$ on $E$. Hence, the collection of all fixed points of $M_\psi$ is preciously the closed subspace $\{f\in \mathbb{T}_p:f|_{E}=0\}$.

\epf

\bcor
Take $E=\{v\in T: \psi(v)\neq 1 \}$. Then, the collection of all fixed points of $M_\psi:\mathbb{T}_{p,0} \rightarrow \mathbb{T}_{p,0}$ is nothing but, the closed subspace $\{f\in \mathbb{T}_{p,0}:f|_{E}=0\}$.
\ecor

\subsection*{Acknowledgments}
The authors thank the National Board for Higher Mathematics (NBHM), India,
for providing financial support  to carry out this research.



\begin{thebibliography} {150}



\bibitem{Colonna-MO-5} R. F. Allen, F. Colonna and G. R. Easley,
\emph{Multiplication operators between Lipschitz-type spaces on a tree},
 Int. J. Math. Math. Sci. (2011), Art. ID 472495, 36 pp.

\bibitem{Colonna-MO-3} R. F. Allen, F. Colonna and G. R. Easley,
\emph{Multiplication operators on the iterated logarithmic Lipschitz spaces of a tree},
 Mediterr. J. Math. {\bf 9} (2012), 575--600.

\bibitem{Colonna-MO-2} R. F. Allen, F. Colonna and G. R. Easley,
\emph{Multiplication operators on the weighted Lipschitz space of a tree},
J. Operator Theory {\bf 69} (2013), 209--231.

\bibitem{Colonna-MO-6} R. F. Allen, F. Colonna and A. Prudhom,
\emph{Multiplication operators between iterated logarithmic Lipschitz spaces of a tree},
Mediterr. J. Math. {\bf 14} (2017) no. 5, Art. 212, 10 pp.


\bibitem{Allen-MO-7} R. F. Allen and I. M. A. Craig,
\emph{Multiplication operators on weighted Banach spaces of a tree},
Bull. Korean Math. Soc. {\bf 54} (3), (2017), 747--761.

\bibitem{Colonna-MO-1} F. Colonna and G. R. Easley,
\emph{Multiplication operators on the Lipschitz space of a tree},
Integr. Equ. Oper. Theory {\bf 68} (2010), 391--411.

\bibitem{Colonna-MO-4} F. Colonna and G. R. Easley,
\emph{Multiplication operators between the Lipschitz space and the space of bounded functions on a tree}, Mediterr. J. Math. {\bf 9} (2012), 423--438.


\bibitem{Conway:Book} J.~B.~ Conway,
{\it A course in functional analysis,} $2$nd edition, Springer, New York, 1990.


\bibitem{MP-Tp-spaces} P. Muthukumar and S. Ponnusamy,
\emph{Discrete analogue of generalized Hardy spaces and multiplication operators on homogenous trees},
Anal. Math. Phys. {\bf 7} (2017), 267--283.


\bibitem{CO-Tp-spaces} P.~Muthukumar and S.~Ponnusamy,
\emph{Composition operators on the discrete analogue of generalized Hardy space on homogenous trees},
Bull. Malays. Math. Sci. Soc. {\bf 40} (4), (2017), 1801--1815.

\bibitem {CO-Tp-spaces2}  P.~Muthukumar and S.~Ponnusamy,
\emph{Composition operators on Hardy spaces of the homogenous rooted trees},
 Monatshefte f$\ddot{u}$r Mathematik (2020), 23 Pages.

\end{thebibliography}
\end{document}